\documentclass[aos,MSNbibl,nameyear,seceqn,dvips]{arximspdf}
\usepackage{dcolumn}

\doi{10.1214/12-AOS983} 
\volume{40}
\issue{3}
\pubyear{2012}
\firstpage{1794}
\lastpage{1815}

\makeatletter

\newcolumntype{d}[1]{D{.}{.}{#1}}

\newproclaim{example}{Example}

\newtheorem{theorem}{Theorem}[section]

\newproclaim{assumption}{Assumption}[section]
\newtheorem{corollary}{Corollary}[section]
\newproclaim{remark}{Remark}[section]

\makeatother

\begin{document}
\begin{frontmatter}

\title{Asymptotic properties of covariate-adaptive randomization}
\runtitle{\hspace*{-10pt}Asymptotic properties of covariate-adaptive randomization}

\begin{aug}
\author{\fnms{Yanqing} \snm{Hu}\ead[label=e1]{yh2s@virginia.edu}}
and
\author{\fnms{Feifang} \snm{Hu}\corref{}\thanksref{t1}\ead[label=e2]{fh6e@virginia.edu}}
\thankstext{t1}{Supported by NSF Grants DMS-09-07297 and DMS-09-06661.}
\runauthor{Y. Hu and F. Hu}
\affiliation{University of Virgina}
\address{Department of Statistics\\
University of Virgina\\
Halsey Hall, Charlottesville\\
Virginia 22904-4135\\
USA \\
\printead{e1}\\
\phantom{E-mail: }\printead*{e2}}
\end{aug}

\received{\smonth{5} \syear{2011}}
\revised{\smonth{2} \syear{2012}}

%
\begin{abstract}
Balancing treatment allocation for influential covariates is critical
in clinical trials.
This has become increasingly important as more and more biomarkers are
found to
be associated with different diseases in translational research
(genomics, proteomics and metabolomics). Stratified permuted block
randomization and
minimization methods [Pocock and Simon \textit{Biometrics} \textbf{31} (1975) 103--115, etc.]
are the two most popular approaches in practice.
However, stratified permuted block randomization
fails to achieve good overall balance when the number of strata is large,
whereas traditional minimization methods also
suffer from the potential drawback of large within-stratum imbalances.
Moreover, the theoretical bases of minimization methods remain largely elusive.
In this paper, we propose a new covariate-adaptive design that is able
to control various types of imbalances.
We show that
the joint process of within-stratum imbalances
is a~positive recurrent Markov chain under certain conditions.
Therefore, this new procedure yields more balanced allocation.
The advantages of the proposed procedure are also demonstrated by
extensive simulation studies.
Our work provides a theoretical tool for future research in this area.
\end{abstract}

%
\begin{keyword}[class=AMS]
\kwd[Primary ]{60F15}
\kwd{62G10}
\kwd[; secondary ]{60F05}
\kwd{60F10}
\end{keyword}

\begin{keyword}
\kwd{Balancing covariates}
\kwd{clinical trial}
\kwd{marginal balance}
\kwd{Markov chain}
\kwd{Pocock and Simon's design}
\kwd{stratified permuted block}
\end{keyword}

\end{frontmatter}

\section{Introduction}
\label{sintro}

In clinical trials, covariates are factors that have a large impact on
the responses of the patients.
Typical covariates include gender, age, disease stage, different
research center, etc.
At the design stage it is often important to balance treatment allocation
over covariates, as a well-balanced trial can lead to more efficient
statistical comparison and more convincing results to the general
audience [\citet{Kun09}]. Balanced allocation is also particularly
useful when the sample size is small or when interim analysis or
subgroup analysis is desired [\citet{Tooetal}].

Stratified randomization is a popular way of achieving balance. It
defines strata as different combinations of the covariates' levels and
employs permuted block randomization within each stratum. This method
is easy to implement and achieves good balance when the number of
strata is small [\citet{KalBeg}]. However, the permuted block
design is susceptible to selection bias [\citet{MatLac88}].
Moreover, it tends to cause severe allocation imbalance in the whole
trial when there are too many strata, typically as a result of many
covariates, or many levels within the individual covariates
[\citet{Poc82}]. Increasing numbers of strata, however, has become the trend,
due to the need to conduct multicenter trials as well as the inclusion
of newly identified biomarkers as covariates [\citet{Khaetal10},
\citet{Lietal10}, \citet{McIetal10}, etc.].

Covariate-adaptive randomization (or minimization) has been
proposed to address the above problem. The earliest work on
minimization dates back to \citet{Tav74} and \citet{PocSim75}.
In particular, with $I$ being the number of covariates and $m_{i}$ the
number of levels for covariate $i$, $i=1,\ldots, I$, Pocock and Simon's
(\citeyear{PocSim75}) procedure minimizes a weighted average of
marginal imbalances $\sum_{i}w_{i}d_{i}(n)$, where $d_{i}(n)$ is a
measure of imbalance among treatment groups with respect to the $i$th
margin of the new patient. Simulation studies [\citet{WeiLee03},
\citet{Tooetal}, \citet{Kun09}] found that this method
reduces marginal imbalances as well as the overall imbalance.
\citet{Wei78} generalized Taves's method by introducing a marginal
urn model. Other works include \citet{Zel74}, \citet
{NorBra77}, Signorini et al. (\citeyear{Sigetal93}) and Heritier, Gebski and
Pillai
(\citeyear{HerGebPil05}), which used a hierarchical decision rule and set priority among
elements of strata, margins and overall trial. Despite the numerous
works in the literature, ``very little is known about the theoretical
properties of covariate-adaptive designs'' [\citet{RosSve08}].

Model-based approach was introduced by \citet{BegIgl80} and
\citet{Atk82}, and the theoretical work has been developed by
Smith (\citeyear{Smi84N1}, \citeyear{Smi84N2}). Smith considered the
linear model $Ey_{n}=T_{n}\alpha+\sum_{j=1}^{p}  z_{n,j}\beta_{j}$
with homogeneous errors and no interaction of any type, where $y_n$,
$T_n$, $(z_{n,1},\ldots,z_{n,p})$ are the response, assignment and
covariate values of the $n$th patient, respectively, and $T_n=+1$ or
$-1$ for treatment 1 or 2. Since $\alpha$, the treatment effect, is the
main interest of the trial, this method sequentially skews the
allocation probability toward the treatment that would lead to a
smaller variance of $\hat{\alpha}$ (the MLE of $\alpha$). Under some
appropriate allocation functions Smith derived the asymptotic normality
of $\sum_{i=1}^{n}z_{i,j}T_{i}$ ($j=1,\ldots, p$). This asymptotic
property was further applied to the construction of a conditional
permutation test [\citet{Smi84N2}].

Although the minimization approach [\citet{PocSim75}, \citet
{Wei78}, etc.] and the model-based approach [\citeauthor{Smi84N1} (\citeyear{Smi84N1,Smi84N2}), etc.]
both lead to marginal and overall balance, they are
rather different in nature. First, even if they use the same biased
coin function, the two allocation rules are still not the same, unless
in the trivial case of no covariates. Hence Smith's asymptotic result
does not readily apply to Pocock and Simon's or Wei's procedure.
Second, Smith's result depends on the homogeneous linear model.
Therefore, once the data type has changed (such as binary or survival
responses), model-based approach does not necessarily imply balance.
Finally, minimization approach is more popular in practice [\citet
{Tav10}]. In fact, as discussed by many authors such as \citet
{LagPoc84}, \citet{Smi84N2} and \citet{McE03}, balanced
allocation enhances credibility of the trials for medical professions
that are less statistically sophisticated, and the simple comparisons
of similar groups of patients are often more acceptable than a
model-based approach adjusting for covariates.

In this paper we focus on the minimization approach that compares patient
numbers at different levels. While the marginal procedures have good
balance with respect to the margins and the whole trial, their
performance within the individual strata is not as satisfactory
[\citet{Sigetal93}, \citet{Kun09}]. \citet{Wei78}
gave a
short proof that if no interaction exists, marginal balances are
sufficient to ensure unbiased estimation of treatment effect in an
unadjusted analysis. In other words, when interactions do exist,
ignorance of within-stratum imbalances may lead to biased estimation.
Moreover, as the field of personalized medicine develops [\citet
{Hu}], subgroup analysis is often desired, and allocation balance
within individual strata can improve the precision of such analysis.

To overcome the potential drawbacks of stratification and Pocock and
Simon's (\citeyear{PocSim75})
method, we develop a new randomization procedure in this paper, which
considers a weighted average of three types of imbalances (within-stratum,
within-covariate-margin and overall). By adopting Efron's (\citeyear
{Efr71}) discrete allocation
function, the next patient will be
assigned with higher probability to a treatment that leads to a smaller
value of the weighted average.

To study the theoretical properties of the new procedure,
the main difficulties include the correlation structure of
within-stratum imbalances as well as the discreteness of the allocation
function. In the literature, a large number of adaptive designs adopt a
continuous allocation function,
and their properties are often obtained by a Taylor expansion
of the allocation function, accompanied by a martingale approximation
[\citet{BaiHu99}, \citet{HuZha04}, \citet{ZhaHuChe06},
etc.]. Since we use Efron's function, which is discrete at 0, the Taylor
expansion is not feasible. We seek to take advantage of an alternative
technique, namely ``drift conditions,'' which was developed to study
the stability of Markov chains on general state spaces. We show that
the joint process of within-stratum imbalances under the new procedure
is a \textit{positive recurrent} Markov chain under some conditions, and
thus preserves the order of $O_{p}(1)$ at the within-stratum level. Our
simulations suggest that the within-stratum imbalances under Pocock and
Simon's (\citeyear{PocSim75}) design have fast-increasing variances as
sample size increases, implying a slower rate than $O_{p}(1)$.

In Section~\ref{spro}, the new procedure is described in general
with $I$ covariates. The theoretical results of the new procedure are
given in
Section~\ref{stheory}. We further use simulations to study the different
covariate-adaptive designs in Section~\ref{ssimu} and conclude our
paper with some observations in Section~\ref{sconclusion}.
The proofs of the theorems can be found in Section~6 and the
supplemental article [\citet{HuHu}].


\section{The new covariate-adaptive randomization procedure}
\label{spro}

This setting is similar to that of \citet{PocSim75}, except that
we only focus on two treatment groups, $1$ and $2$. Consider $I$
covariates and $m_{i}$ levels
for the $i$th covariate, resulting in $m=\prod_{i=1}^{I}m_{i}$ strata.
Let $T_{j}$ be the assignment of the $j$th patient, $j=1,\ldots,n$,
that is, $T_{j}=1$ for treatment 1 and $T_{j}=0$ for treatment 2.
Let $Z_{j}$ indicate the covariate profile of that patient,
that is, $Z_{j}=(k_{1},\ldots,k_{I})$ if his or her $i$th covariate is
at level $k_{i}$, $1\leq i \leq I $ and $1\leq k_{i}\leq m_{i}$.
For convenience, we use $(k_{1},\ldots,k_{I})$ to denote the \textit
{stratum} formed by patients who possess the same covariate profile
$(k_{1},\ldots,k_{I})$, and use $(i;k_{i})$ to denote the \textit{margin}
formed by patients whose $i$th covariate is at level $k_{i}$.

The new procedure is defined as follows:
\begin{enumerate}[(1)]
\item[(1)] The first patient is assigned to treatment 1 with
probability 1/2.

\item[(2)] Suppose $(n-1)$ patients have been assigned to a treatment
($n>1$) and
the $n$th patient falls within stratum $(k_{1}^{*},\ldots,k_{I}^{*})$.

\item[(3)] For the first $(n-1)$ patients:
\begin{enumerate}
\item[-]let $D_{n-1}$ be the difference between the numbers of patients
in treatment group 1 and 2 as total, that is, the number in group 1
minus the number in group 2;
\item[-]similarly, let $D_{n-1}(i;k_{i}^{*})$ and
$D_{n-1}(k_{1}^{*},\ldots,k_{I}^{*})$ be the differences between the
numbers of patients in the two treatment groups on the margin
$(i;k_{i}^{*})$, and within the stratum $(k_{1}^{*},\ldots,
k_{I}^{*})$, respectively;
\item[-] these differences can be positive, negative or zero, and each
one is used to measure the \textit{imbalance} at the corresponding level
(overall, marginal, or within-stratum).
\end{enumerate}

\item[(4)] If the $n$th patient were assigned to treatment $1$, then
$D^{(1)}_{n}=D_{n-1}+1$ would
be the ``potential'' overall difference in the two groups; similarly,
\[
D^{(1)}_{n}\bigl(i;k_{i}^{*}\bigr)=
D_{n-1}\bigl(i;k_{i}^{*}\bigr)+1
\]
and
\[
D^{(1)}_{n}\bigl(k_{1}^{*},
\ldots,k_{I}^{*}\bigr)=D_{n-1}\bigl(k_{1}^{*},
\ldots,k_{I}^{*}\bigr)+1
\]
would be the potential differences on margin $(i;k_{i}^{*})$ and within
stratum $(k_{1}^{*},\ldots,k_{I}^{*})$, respectively.\vadjust{\goodbreak}

\item[(5)] Define an imbalance measure $\mathit{Imb}_{n}^{(1)}$ by
\[
\mathit{Imb}_{n}^{(1)}=w_{o}
\bigl[D^{(1)}_{n}\bigr]^2+\sum
_{i=1}^{I}w_{m,i}\bigl[D^{(1)}_{n}
\bigl(i;k_{i}^{*}\bigr)\bigr]^2 +w_{s}
\bigl[D^{(1)}_{n}\bigl(k_{1}^{*},
\ldots,k_{I}^{*}\bigr)\bigr]^2,
\]
which is the weighted imbalance that would be caused if the $n$th patient
were assigned to treatment 1. $w_{o}$, $w_{m,i}$ and $w_{s}$ are nonnegative
weights placed overall, within a covariate margin and within a stratum
cell, respectively. Without loss of generality we can assume
\[
w_{o}+w_{s}+\sum_{i=1}^{I}w_{m,i}=1.
\]
\item[(6)] In the same manner we can define $\mathit{Imb}_{n}^{(2)}$, the
weighted imbalance that would
be caused if the $n$th patient were assigned to treatment 2. In this
case, the three types
of potential differences are the existing ones minus 1, instead of plus 1.

\item[(7)] Conditional on the assignments of the first $(n-1)$ patients
as well as the covariates' profiles of the first $n$ patients, assign
the $n$th patient to treatment 1 with probability
\[
P(T_{n}=1|\mathbf{Z}_{n},\mathbf{T}_{n-1})=
\cases{ %
q, & \quad $\mbox{if $\mathit{Imb}_{n}^{(1)}>
\mathit{Imb}_{n}^{(2)}$,}$
\vspace*{2pt}\cr
p, & \quad$\mbox{if $\mathit{Imb}_{n}^{(1)}<
\mathit{Imb}_{n}^{(2)},$}$
\vspace*{2pt}\cr
0.5, & \quad$\mbox{otherwise,}$}
\]
where $n>1$, $0<q<p<1$, $p+q=1$, $\mathbf{Z}_{n}=(Z_{1},\ldots,Z_{n})$
and $\mathbf{T}_{n-1}=(T_{1},\ldots,T_{n-1})$.
\end{enumerate}
%
%
\begin{remark}
When $w_{o}=w_{s}=0$, that is, only the marginal imbalances are
considered, the proposed design
reduces to a special case of Pocock and Simon's (\citeyear{PocSim75})
method; and when $w_{m,i}=w_{o}=0$,
it reduces to stratified randomization, where a separate biased coin is
employed to determine
the assignment within each stratum. However, we will explore procedures
with other choices of weights, to see if they can lead to more balanced
allocation from various perspectives.
\end{remark}

%
\begin{remark}
In the literature different views have been given as to the selection
of the biasing probability $p$. \citet{Efr71} suggested $p=2/3$,
but his method does not consider covariates. The more recent papers,
especially those involving covariate-adaptive randomization, suggested
larger $p$'s, such as 0.85, 0.90 and 0.95. See \citet{WeiLee03},
Hagino et al. (\citeyear{Hagetal04}), Toorawa et al. (\citeyear{Tooetal}), and
\citet{HuZhaHe09}. One may also use other generators in step (7),
for example, Wei's (\citeyear{Wei78}) generator.
The properties of the design will be different.
\end{remark}

%
\begin{example}
Suppose in a trial two covariates, gender and smoking behavior, are
considered to be influential,
each of which has two levels. Thus, the 4 strata $(1,1)$, $(1,2)$,
$(2,1)$, $(2,2)$ represent male smokers, male nonsmokers, female
smokers and female nonsmokers, respectively. Assume that the weights are
$w_{o}=1/3$, $w_{m,1}=w_{m,2}=1/6$ and $w_{s}=1/3$. The first 50
patients have been randomized
and the 4 within-stratum differences among these 50 patients are $-2$,
$+2$, $+1$ and $-1$. If the 51th patient is a \textit{male smoker}, then
the current imbalances are calculated as:
\begin{itemize}
\item[-] overall: $D_{n-1}=-2+2+1-1=0$;
\item[-] margin of male: $D_{n-1}(1;1)=-2+2=0$;
\item[-] margin of smokers: $D_{n-1}(2;1)=-2+1=-1$;
\item[-] stratum of male smokers: $D_{n-1}(1,1)=-2$.
\end{itemize}
The potential imbalances if the new patient were assigned to treatment
1 or~2 are given in Table~\ref{texample50patients}.
%
%
\begin{table}
\caption{An example illustrating the calculation under the new procedure}
\label{texample50patients}
\begin{tabular*}{\textwidth}{@{\extracolsep{\fill}}lccc@{}}
\hline
&$\bolds{D_{50}(\cdot)}$ &$\bolds{D_{51}^{(1)}(\cdot)=D_{50}(\cdot)+}1$
&$\bolds{D_{51}^{(2)}(\cdot)=D_{50}(\cdot)-1}$\\
\hline
Overall &\phantom{0}0 &\phantom{0}1 &$-1$\\
Margin of male (1;1) &\phantom{0}0 &\phantom{0}1 &$-1$\\
Margin of smokers (2;1) &$-1$ &\phantom{0}0 &$-2$\\
Stratum of male smoker (1,1) &$-2$ &$-1$ &$-3$\\
\hline
\end{tabular*}
\end{table}

Therefore,
\begin{eqnarray*}
\mathit{Imb}_{51}^{(1)}&=&(1)^2\cdot
\tfrac{1}{3}+(1)^2 \cdot\tfrac
{1}{6}+(0)^2\cdot
\tfrac{1}{6}+(-1)^2 \cdot\tfrac{1}{3}=0.83,
\\
\mathit{Imb}_{51}^{(2)}&=&(-1)^2\cdot
\tfrac{1}{3}+(-1)^2 \cdot\tfrac
{1}{6}+(-2)^2
\cdot\tfrac{1}{6}+(-3)^2 \cdot\tfrac{1}{3}=4.17
.
\end{eqnarray*}
Since $\mathit{Imb}_{51}^{(1)}=0.83<\mathit{Imb}_{51}^{(2)}=4.17$, the coin will be
biased toward treatment 1 with probability $p>0.5$.
\end{example}

\section{Theoretical properties of the new design}
\label{stheory}

We now investigate the asymptotic properties of the proposed design.
For the first $n$ patients, we know that $D_{n}(k_{1},\ldots,k_{I})$
is the true difference of patient numbers within stratum $(k_{1},\ldots
,k_{I})$. Furthermore, let
\[
\mathbf{D}_{n}= \bigl[D_{n}(k_{1},
\ldots,k_{I}) \bigr]_{1\leq
k_{1}\leq
m_{1},\ldots, 1\leq k_{I}\leq m_{I}}
\]
be an array of dimension $m_{1}\times\cdots\times m_{I}$ which stores
the current assignment
differences in all strata.
Also, assume that the covariates $Z_{1},Z_{2},\ldots$ are independently
and identically
distributed. Since $Z_{n}=(k_{1},\ldots, k_{I})$ can take
$m=\prod_{i=1}^{I} m_i$ different values, it in fact follows an $m$-dimension
multinomial distribution\vadjust{\goodbreak} with parameter $\mathbf{p}=(p(k_{1},\ldots,
k_{I}))$, each element being the probability that a patient falls
within the corresponding stratum. Obviously, $p(k_{1},\ldots, k_{I})
\geq0 $ and $\sum_{k_{1},\ldots, k_{I}} p(k_{1},\ldots, k_{I})=1$.\vspace*{1pt}

First, we notice that $(\mathbf{D}_{n})_{n\geq1}$ is a Markov chain on
the space $\mathbb{Z}^{m}$. In fact, by definition of the new
procedure, $\mathbf{D}_{n}$ is a function $f$ of $(\mathbf
{D}_{n-1},Z_{n}, T_{n})$. Moreover, conditional on $\mathbf{D}_{n-1}$,
$(Z_{n},T_{n})$ is independent of $(\mathbf{D}_{1},\ldots,\mathbf
{D}_{n-2})$; therefore, $\mathbf{D}_{n}=f(\mathbf{D}_{n-1},Z_{n},
T_{n})$ is also conditionally independent of $(\mathbf{D}_{1},\ldots,\break\mathbf{D}_{n-2})$.

We next explore the conditions under which $(\mathbf{D}_{n})_{n\geq1}$
is a \textit{positive recurrent} chain, a desired property which
indicates fast convergence rate. We will first investigate the special
case of $2\times2$ strata, that is, only two covariates and two levels
for each. This case enables us to obtain a finer result than the more
general case, and at the same time also sheds light on how to set the
conditions for the latter. With $2\times2$ strata, the weights on
$\mathit{Imb}^{(1)}_{n}$ or $\mathit{Imb}^{(2)}_{n}$ reduce to $w_{o}$, $w_{m,1}$,
$w_{m,2}$ and~$w_{s}$.

%
\begin{theorem}\label{theorem1}
For the new design, consider 2 covariates and 2 levels for each.
$w_{o}$, $w_{m,1}$, $w_{m,2}$ and $w_{s}$ are nonnegative with
$w_{o}+w_{m,1}+w_{m,2}+w_{s}=1$. If the following two conditions hold:
\begin{longlist}[(A)]
\item[(A)]$w_{s}>0$,

\item[(B)] define
\begin{eqnarray*}
u_{1}&=&w_{o}+w_{m,1}+w_{m,2}+w_{s}=1,
\\[-1pt]
u_{2}&=&w_{o}+w_{m,1},
\\[-1pt]
u_{3}&=&w_{o}+w_{m,2},
\\[-1pt]
u_{4}&=&w_{o};
\end{eqnarray*}
the solution $\mathbf{x}=(x_{1},x_{2},x_{3})$ to the linear equation
\[
\pmatrix{ u_{1} & u_{2} &
u_{3}
\vspace*{2pt}\cr
u_{2} & u_{1} & u_{4}
\vspace*{2pt}\cr
u_{3} & u_{4} & u_{1}
}
\pmatrix{
x_{1}
\vspace*{2pt}\cr
x_{2}
\vspace*{2pt}\cr
x_{3}
}
= \pmatrix{
u_{4}
\vspace*{2pt}\cr
u_{3}
\vspace*{2pt}\cr
u_{2}
}
\]
satisfies $|x_{1}|+|x_{2}|+|x_{3}|<1$,
\end{longlist}
then $(\mathbf{D}_{n})_{n\geq1}$ is a positive recurrent Markov chain
with period 2 on $\mathbb{Z}^{4}$.
\end{theorem}

%
\begin{remark}
By Theorem~\ref{theorem1} the chains $\mathbf{D}_{2n+1}$ and $\mathbf
{D}_{2n}$ are two ergodic chains
and converge to two limit distributions, respectively. Thus, $\mathbf
{D}_{2n+1}=O_{p}(1)$ and
$\mathbf{D}_{2n}=O_{p}(1)$, which implies that $\mathbf
{D}_{n}=O_{p}(1)$. Accordingly, the imbalances
at any level (within strata, on the margins, or overall) preserve the
order of $O_{p}(1)$.
\end{remark}
%
%
\begin{remark}
In fact, $u_{1}$, $u_{2}$,
$u_{3}$ and $u_{4}$ in the above theorem can be interpreted as the
weights placed on individual strata: call the stratum in which the\vadjust{\goodbreak}
current patient falls a ``target,'' then $u_{1}$ is the weight on the
target itself; $u_{2}$ ($u_{3}$) on any stratum that is at the same
level of covariate 1 (covariate 2) as the target; $u_{4}$ on any of the
remaining strata.
\end{remark}
%
%
\begin{corollary}\label{corollary1}
In Theorem~\ref{theorem1}, if we further assume that $w_{m,1}=w_{m,2}:=w_{m}$,
then condition \textup{(B)} is equivalent to
%
%
{\renewcommand{\theequation}{$\mathrm{B}^{\prime}$}
\begin{equation}\label{bprime}
w_{m}<C(w_{o}):=\frac{\sqrt{(1-w_{o})^{2}+4(1+w_{o})^{2}}-1-3w_{o}}{4}.\vspace*{-3pt}
\end{equation}
}
\end{corollary}

%
\begin{table}
\caption{Constraint on $w_{m}$ as a function of $w_{o}$}
\label{tconstraint}
\begin{tabular*}{200pt}{@{\extracolsep{\fill}}lccccc@{}}
\hline
$w_{o}$ &0.00 &0.20 &0.40 &0.60 &0.80 \\
$C(w_{o})$ &0.31 &0.23 &0.17 &0.11 &0.05\\
\hline
\end{tabular*}     \vspace*{-3pt}
\end{table}

%
%
\begin{remark}
See Table~\ref{tconstraint} for certain values of $C(w_{o})$. Since
$C(w_{o})$ is a~decreasing and almost linear function of $w_{o}$ on
$[0,1]$, condition (\ref{bprime})
is much easier to verify than condition (B). For example, if the
weight at the overall level $w_{o}=0.20$, then the ones on the two
margins need to be less than 0.23. Therefore,
$(w_{o},w_{m,1},w_{m,2},w_{s})$ = $(0.20,0.22,0.22,0.36)$ is a
legitimate weight set that ensures positive recurrence.\vspace*{-3pt}
\end{remark}

The next theorem deals with the general case of $m=\prod_{i=1}^{I}m_{i}$ strata. Using the basic equation $(x+1)^2-(x-1)^2=4x$,
the critical quantity $\mathit{Imb}_{n}^{(1)}-\mathit{Imb}_{n}^{(2)}$ in step $(7)$
(Section~\ref{spro}) can be simplified as
%
%
\renewcommand{\theequation}{\arabic{section}.\arabic{equation}}
\setcounter{equation}{0}
\begin{eqnarray}\label{simplification}
&&\mathit{Imb}_{n}^{(1)}-\mathit{Imb}_{n}^{(2)}
\nonumber\\
&&\hspace*{2pt}\qquad= 4 \Biggl\{w_{o}D_{n-1}+\sum_{i=1}^{I}w_{m,i}D_{n-1}
\bigl(i;k_{i}^{*} \bigr)+w_{s}D_{n-1}
\bigl(k_{1}^{*},\ldots,k_{I}^{*} \bigr)
\Biggr\}
\\
&&\qquad:= 4\cdot\delta_{n-1} \bigl(k_1^*,\ldots,k_I^*
\bigr) .\nonumber
\end{eqnarray}
Therefore, the biasing probability $p$, $q$ or $1/2$ is determined by
the sign of $\delta_{n-1}(k_1^*,\ldots,k_I^*)$, which is a weighted
average of current imbalances at different levels. Since $D_{n-1}$ and
$D_{n-1}(i;k_{i}^{*})$ can both be expressed as a~sum of certain
$D_{n-1}(k_{1},\ldots,k_{I})$'s, we want to reformulate
$\delta_{n-1}(k_1^*,\ldots,k_I^*)$ as a~weighted average of imbalances within
the individual strata.

As a motivating example, consider 3 covariates, gender (male or
female), smoking behavior (smoker or nonsmoker) and clinical center (3
centers), with a total of 12 strata. Suppose for the new patient
$Z_{n}=(1,1,1)$, that is, he falls into the stratum of ``male smokers
at center 1.'' Then for the remaining strata, the weights on
$D_{n-1}(k_{1},\ldots,k_{I})$'s in the expression of $\delta_{n-1}(k_1^*,\ldots,k_I^*)$
are shown in Table~\ref{texampleweights}.

%
%
\begin{table}
\caption{An example showing the weights of $D_{n-1}(k_{1},\ldots
,k_{I})$'s in $\delta_{n-1}(k_1^*,\ldots,k_I^*)$}
\label{texampleweights}
\begin{tabular*}{\textwidth}{@{\extracolsep{\fill}}lccc@{}}
\hline
&\textbf{Stratum} &\textbf{Description} &\textbf{Weight}\\
\hline
\phantom{0}1 & (1,1,1) & male smokers at center 1 &
$w_{o}+w_{m,1}+w_{m,2}+w_{m,3}+w_{s}=1$\\
\phantom{0}2 & (1,1,2) & male smokers at center 2 &
$w_{o}+w_{m,1}+w_{m,2}$\\
\phantom{0}3 & (1,1,3) & male smokers at center 3 &
$w_{o}+w_{m,1}+w_{m,2}$\\
\phantom{0}4 & (1,2,1) & male nonsmokers at center 1 &
$w_{o}+w_{m,1}+w_{m,3}$\\
\phantom{0}5 & (2,1,1) & female smokers at center 1 &$
w_{o}+w_{m,2}+w_{m,3}$\\
\phantom{0}6 & (1,2,2) & male nonsmokers at center 2 &
$w_{o}+w_{m,1}$\\
\phantom{0}7 & (1,2,3) & male nonsmokers at center 3 &
$w_{o}+w_{m,1}$\\
\phantom{0}8 & (2,1,2) & female smokers at center 2 & $w_{o}+w_{m,2}$\\
\phantom{0}9 & (2,1,3) & female smokers at center 3 & $w_{o}+w_{m,2}$\\
10& (2,2,1) & female nonsmokers at center 1 & $w_{o}+w_{m,3}$\\
11& (2,2,2) & female nonsmokers at center 2 & $w_{o}$\\
12& (2,2,3) & female nonsmokers at center 3 & $w_{o}$\\
\hline
\end{tabular*}
\end{table}

Generally, with respect to stratum $(k_{1}^{*},\ldots,k_{I}^{*})$ in
which the new patient falls, we will divide\vadjust{\goodbreak} the $m=\prod_{i=1}^{i}m_{i}$ strata into several categories and find out the
corresponding weights in the expression of $\delta_{n-1}(k_1^*,\ldots
,k_I^*)$. Let $\mathbb{I}=\{1,2,\ldots, I\}$. For any stratum
$(k_{1},\ldots,k_{I})$:
\begin{itemize}
\item[-] if $(k_{1},\ldots,k_{I})=(k_{1}^{*},\ldots,k_{I}^{*})$, then
the weight on $D_{n-1}(k_{1},\ldots,k_{I})$ is $w_{o}+\sum_{i=1}^{I}w_{m,i}+w_{s}=1$;

\item[-] for any fixed $i$ ($i\in\mathbb{I}$), if $k_{i}\neq
k_{i}^{*}$ and $k_{j}=k_{j}^{*}$ for $j\in\mathbb{I}$ and $j\neq i$,
then the weight on $D_{n-1}(k_{1},\ldots,k_{I})$ is $w_{o}+\sum_{j\neq
i}w_{m,j}$, and there are $(m_{i}-1)$ strata in this category;

\item[-] for any fixed $i_{1}< i_{2}$ ($\{i_{1},i_{2}\}\subset\mathbb
{I}$), if $k_{i_{1}}\neq k_{i_{1}}^{*}$, $k_{i_{2}}\neq k_{i_{2}}^{*}$,
and $k_{j}=k_{j}^{*}$ for $j\in\mathbb{I}$, $j\neq i_{1}$ and $j\neq
i_{2}$, then the weight on $D_{n-1}(k_{1},\ldots,k_{I})$ is
$w_{o}+\sum_{j\neq i_{1} , j\neq i_{2}}w_{m,j}$, and there are
$(m_{i_{1}}-1)(m_{i_{2}}-1)$ strata in this category;

\item[-] for any fixed $i_{1}< i_{2}< \cdots< i_{l}$ ($\{i_{1},\ldots
,i_{l}\}\subset\mathbb{I}$), if $k_{i_{t}}\neq k_{i_{t}}^{*}$ and
$k_{j}=k_{j}^{*}$ for $j\in\mathbb{I}$, $j\neq i_{t}$ and $1\leq t
\leq l$, then the weight on $D_{n-1}(k_{1},\ldots,k_{I})$ is
$w_{o}+\sum_{j\neq i_{t}, 1\leq t \leq l}w_{m,j}$, and there are
$\prod_{t=1}^{l}(m_{i_{t}}-1)$ strata in this category;

\item[-] if $k_{i}\neq k_{i}^{*}$ for all $i\in\mathbb{I}$, then the
weight on $D_{n-1}(k_{1},\ldots,k_{I})$ is $w_{o}$, and there are
$\prod_{i=1}^{l}(m_{i}-1)$ strata in this category.
\end{itemize}
It is easily verified that
\begin{eqnarray*}
m&=&\prod_{i=1}^{I}m_{i}
\\
&=& \bigl[(m_{1}-1)+1 \bigr] \bigl[(m_{2}-1)+1 \bigr]\cdots
\bigl[(m_{I}-1)+1 \bigr]
\\
&=&1+\sum_{l=1}^{I}\sum
_{1\leq i_{1}< i_{2}<\cdots< i_{l}\leq I}\prod_{t=1}^{l}[m_{i_{t}}-1],
\end{eqnarray*}
which is consistent with the counts listed above. Our general theorem
in the following is closely related to the above weights and counts.
%
%
\begin{theorem}\label{theorem2}
For the new design, consider $I$ covariates and $m_{i}$ levels for the
$i$th covariate, where $I\geq1$, $1\leq i \leq I$, and $m_{i}>1$.
$w_{o}$, $w_{s}$ and $w_{m,i}$, $i=1,\ldots,I$, are nonnegative with
$w_{o}+\sum_{i=1}^{I}w_{m,i}+w_{s}=1$. If
%
%
{\renewcommand{\theequation}{$\mathrm{C}$}
\begin{equation}\label{eqC}
u^{*}:=\sum_{l=1}^{I}\sum
_{1\leq i_{1}< i_{2}<\cdots< i_{l}\leq
I} \Biggl\{ \biggl(w_{o}+\sum
_{j\neq i_{t},
1\leq t \leq l}w_{m,j} \biggr)\prod
_{t=1}^{l}[m_{i_{t}}-1] \Biggr\}<1/2,\hspace*{-35pt}
\end{equation}}
\hspace*{-2pt}then $\mathbf{D}_{n}$ is a positive recurrent Markov chain on $\mathbb
{Z}^{m}$.
\end{theorem}

To see the theorem in a more intuitive way, we will take a closer look
at $u^{*}$ in the special case of two covariates, as is shown in the
following corollary.
%
%
\begin{corollary}\label{corollary2}
In Theorem~\ref{theorem2}, if $I=2$, then condition (\textup{\ref{eqC}}) is equivalent to
%
%
{\renewcommand{\theequation}{$\mathrm{C}^{\prime}$}
\begin{equation}\label{eqCprime}
(m_{1}m_{2}-1)w_{o}+(m_{1}-1)w_{m,2}+(m_{2}-1)w_{m,1}<1/2.
\end{equation}}
\end{corollary}
%
%
\begin{remark}
When $w_{o}=0$ and $w_{m,1}=w_{m,2}=w_{m}$, condition \textup{(\ref{eqCprime})} further
reduces to $w_{m}<[2(m_{1}+m_{2}-2)]^{-1}$. For example, if
$m_{1}=m_{2}=5$, then $w_{m}<1/16$ is required to satisfy condition (\textup{\ref{eqC}}).
\end{remark}
%
%
\begin{remark}
In both Theorems~\ref{theorem1} and~\ref{theorem2}, $w_{s}>0$
is required. Therefore, the theoretical results in these theorems do
not apply to Pocock and Simon's (\citeyear{PocSim75}) design (with
$w_s=0$). The simulation result in Table~\ref{t2times2table}
(Section~\ref{ssimu}) shows that the within-stratum imbalances under their
method increase as the sample size increases, suggesting that they may
not have the rate of $O_{p}(1)$. We hypothesize that the condition
$w_s>0$ is critical to ensure that $(\mathbf{D}_{n})_{n\geq1}$ is
positive recurrent. These are further research problems.
\end{remark}

To prove the above two theorems, we will use the technique of ``drift
conditions'' [\citet{MeyTwe93}], which was developed for Markov
chains on general state spaces. Applying their theory to our problem,
in order to prove positive recurrence of $(\mathbf{D}_{n})_{n\geq1}$
we need to find a test function
$V\dvtx\mathbb{Z}^{m}\rightarrow\mathbb{R}^{+}$, a bounded test set
$\mathcal{C}$ on $\mathbb{Z}^{m}$
and two positive constants $M_{1}$ and $M_{2}$ such that
%
%
\renewcommand{\theequation}{\arabic{section}.\arabic{equation}}
\setcounter{equation}{1}
\begin{equation}
\Delta V (\mathbf{D}):= \sum_{\mathbf{D}'\in\mathbb
{Z}^{m}}P \bigl(\mathbf
{D},\mathbf{D}' \bigr)V \bigl(\mathbf{D}' \bigr)-V(
\mathbf{D})\label{deltavd}
\end{equation}
satisfies the following two conditions:
%
%
\begin{eqnarray}
\Delta V (\mathbf{D})&\leq&-M_{1},\qquad \mathbf{D}\notin\mathcal{C},
\label{drift1}
\\
\Delta V (\mathbf{D})&\leq& M_{2},\qquad \mathbf{D}\in\mathcal{C},
\label{drift2}
\end{eqnarray}
where $P(\mathbf{D},\mathbf{D}')$ is the transition probability from
$\mathbf{D}$ to $\mathbf{D}'$
on state space~$\mathbb{Z}^{m}$ of the chain $(\mathbf{D}_{n})_{n\geq
1}$. $V$ is often a norm-like function
on $\mathbb{Z}^{m}$. These drift conditions can roughly be interpreted
as follows: so long as the average
one-step movement $\Delta V$ tends to go back (with the magnitude
uniformly greater than a positive
constant~$M_{1}$), that is, the chain is pulled back toward the finite
set $\mathcal{C}$, positive recurrence can be ensured. For proofs of
the theorems, see Section~\ref{ssketch} and the supplemental article [\citet{HuHu}].

\section{Simulation studies}
\label{ssimu}

We will compare the new procedure with stratified permuted block
randomization and Pocock and Simon's (\citeyear{PocSim75}) minimization
method. The simulations can be divided into three parts. First, we will
simulate the case of $2\times2$ strata with a relatively large number
of patients, to verify the convergence rate as stated in Theorem~\ref
{theorem1}. Secondly, we are interested in the performances of
different randomization methods when the number of strata is large as
compared to the sample size. An example of 500 patients and 10
covariates (each with 2 levels) will be studied. Finally, an example
from Toorawa et al. (\citeyear{Tooetal}) will be considered, which is chosen
because it resembles real situations in clinical trials.

\subsection{\texorpdfstring{$2\times2$ strata}{2x2 strata}}
For the three randomization procedures, we want to see whether the
imbalances at any of the three levels (within-stratum, marginal and
overall) stabilize, which indicates the rate of $O_{p}(1)$ at that
specific level. The parameters are specified as follows:
\begin{itemize}
\item[-] Multinomial probability
$(p(1,1),p(1,2),p(2,1),p(2,2))=(0.1,0.2,0.3,0.4)$.
\item[-] Biasing probability $p=0.85$ and $q=0.15$ for Pocock and
Simon's method (PS) as well as for the new procedure (NEW).
\item[-] Block size 4 for stratified randomization (STR-PB).
\item[-] Sample size $n=200,500,1000$; number of simulated trials $N=1000$.
\item[-] NEW: $(w_{o},w_{m,1},w_{m,2},w_{s})=(0.3,0.1,0.1,0.5)$;
conditions \textup{(A)} and \textup{(B)} are satisfied.
\item[-] PS: $(w_{o},w_{m,1},w_{m,2},w_{s})=(0,0.5,0.5,0)$; Conditions
\textup{(A)} and \textup{(B)} are NOT satisfied.
\end{itemize}
Table~\ref{t2times2table} shows the standard deviations (std's) of
$D_{n}(\cdot)$'s at different levels (by symmetry of the designs, the
theoretical
mean of each $D_{n}(\cdot)$ is always~$0$). For simplicity, only the
result of 2 strata and 2 margins are listed. Of the five columns, the
first and the second give the std's of assignment differences within
stratum $(1,1)$ and $(2,2)$; the third and fourth for the marginal
differences of covariate 1 at level 1 and covariate 2 at level 2; and
the last for the overall difference.

%
\begin{table}
\caption{std's of $D_n(\cdot)$ of several methods under different
sample sizes}
\label{t2times2table}
\begin{tabular*}{\textwidth}{@{\extracolsep{\fill}}lcccccc@{}}
\hline
\multicolumn{2}{@{}l}{\textbf{Sample size}}&$\bolds
{D_{n}(1,1)}$&$\bolds{D_{n}(2,2)}$&$\bolds{D_{n}(1;1)}$&$\bolds
{D_{n}(2;2)}$&\multicolumn{1}{c@{}}{$\bolds{D_{n}}$}\\
\hline
STR-PB&$\phantom{0}200$&0.92&0.89&1.30&1.27&1.83 \\
&\phantom{0}$500$&0.92&0.92&1.31&1.30&1.86 \\
&$1000$&0.92&0.89&1.31&1.28&1.81\\[3pt]
PS&\phantom{0}$200$&3.16&3.27&1.15&1.13&1.30\\
&\phantom{0}$500$&4.80&4.83&1.16&1.11&1.31\\
&$1000$&7.25&7.33&1.15&1.13&1.30\\[3pt]
NEW&\phantom{0}$200$&1.11&1.07&1.30&1.27&1.32\\
&\phantom{0}$500$&1.14&1.10&1.33&1.28&1.22\\
&$1000$&1.03&1.10&1.20&1.24&1.27 \\
\hline
\end{tabular*}
\end{table}

Table~\ref{t2times2table} suggests that all 5 standard deviations
stabilize under NEW and \mbox{STR-PB} when the sample size increases. For
example, under NEW the std's of $D_{n}(1,1)$
are~1.11, 1.14 and 1.03; and those of $D_{n}$ are 1.32, 1.22 and 1.27,\vadjust{\goodbreak}
which means that our new
procedure preserves the rate of $O_{p}(1)$. The same conclusion can be
reached for STR-PB. In fact,
since the block size is~4, any within-stratum imbalance under STR-PB is
bounded by 2.
For PS, however, while the std's of marginal and overall differences
stabilize, those
of the within-stratum differences do not. For example, the std of
$D_{n}(1,1)$ increases from
3.16 to 4.80 and 7.25, much larger than those under the other two methods.

For the within-stratum imbalances, STR-PB is the best [0.92 for
$D_{n}(1,1)$], with NEW having slightly larger std's and PS the
largest. For the marginal imbalances, PS is the best [around 1.15 for
$D_{n}(1;1)$],
and the other two are about the same [around 1.30 for $D_{n}(1;1)$].
For the overall imbalance,
STR-PB is not as good as NEW and PS. Therefore, we see that even for 4
strata, STR-PB does not perform well
for the overall imbalance.

\subsection{\texorpdfstring{$2^{10}$ strata}{2 10 strata}}

We simulate a hypothetical trial, which involves 500 patients, 10
covariates and 2 levels for each, that is, 1024 strata in total. The
biased coin probabilities $p$ and $q$ for NEW and PS, the block size
for STR-PB and the number of simulated trials $N$ remain the same. The
covariates are generated as follows: in addition to the independence
assumption of covariates between patients, we further assume that
within each patient the different covariates are independent and that
each level within a fixed covariate is equally likely. Therefore, for
the covariate profile $Z_{i}=(k_{1},\ldots, k_{I})$ of the $i$th
patient, $k_{1},\ldots, k_{I}$ are independently sampled from $\{1,2\}$.
For the weights, we use $w_{o}=0$, $w_{s}=0.5$ and $w_{m,i}=0.5/10$.

Of the 1024 strata, on average 61.4\% have no patient, and only 0.1\%
have 4 or more. Hence, if STR-PB is employed, most blocks are
incomplete, which tends to cause large overall imbalance. Table~\ref
{t10241} displays the mean absolute imbalances under each of the
three randomization methods.
%
%
\begin{table}
\caption{Mean $|D_{n}(\cdot)|$ for $2^{10}$ strata and 500 patients}
\label{t10241}
\begin{tabular*}{200pt}{@{\extracolsep{\fill}}lccc@{}}
\hline
&\textbf{STR-PB} &\textbf{PS} & \multicolumn{1}{c@{}}{\textbf
{NEW}}\\
\hline
Overall &17.07 &0.76 &0.98 \\
Marginal &11.80 &1.65 &1.94 \\
Within-strt. (2 pts) &\phantom{0}0.66 &0.98 &0.50 \\
Within-strt. (3 pts) &\phantom{0}1.00 &1.23 &1.08 \\
\hline
\end{tabular*}
\end{table}

As seen in Table~\ref{t10241}, STR-PB has an extremely large
E$|D_{n}|$ (17.07). In comparison,\vadjust{\goodbreak} the other two methods have much
smaller values of 0.76 and 0.98. So in this respect, PS has the best
performance, and NEW is only slightly worse. In the second row, the
mean absolute marginal imbalance is the average of the absolute
differences over 20 margins as well as over the 1000 simulations, and
the interpretation is the same as the overall imbalance. For the
within-stratum imbalances, the table shows the result for strata with 2
or 3 patients. For example, under PS, 0.98 is the mean absolute
difference over all strata with 2 patients as well as over the 1000
simulations. Under this criterion, PS is not recommended since the two
means are 0.98 and 1.23, the largest among the three methods. STR-PB
and NEW are quite similar, with means 0.66 versus 0.50 for strata with~2 patients, and 1.00 versus 1.08 with 3 patients. Hence, although our
new procedure is not always the best, it ensures that no single type of
the imbalances becomes too extreme.
%
%
\begin{table}[b]
\caption{Distribution of covariates}
\label{tDistributionofcovariates}
\begin{tabular*}{\textwidth}{@{\extracolsep{\fill}}lcc@{}}
\hline
Sites & Small (2 sites) &$1/120$\\
& Medium (16 sites) &$6/120$\\
& Large (2 sites) &$11/120$\\
Other 3 covariates &Male; $<60$; Moderate disease & $10/20$\\
&Male; $\geq60$; Moderate disease & $2/20$\\
&Male; $<60$; Severe disease & $2/20$\\
&Male; $\geq60$; Severe disease & $2/20$\\
&Female; $<60$; Moderate disease & $1/20$\\
&Female; $\geq60$; Moderate disease & $1/20$\\
&Female; $<60$; Severe disease & $1/20$\\
&Female; $\geq60$; Severe disease & $1/20$\\
\hline
\end{tabular*}
\end{table}

\subsection{An example mimicking real clinical data}
We chose an example from \citet{Tooetal}. The four covariates are
site, gender, age and disease status, with 20, 2, 2 and 2 levels,
respectively, resulting in 160 strata. The covariates' distribution is
replicated in Table~\ref{tDistributionofcovariates}, where the
marginal distribution of sites is independent of the joint distribution
of the remaining three covariates.

120 patients enter the trial sequentially, and their covariates are
independently simulated
from the multinomial distribution in Table~\ref
{tDistributionofcovariates}. We use the same $p$, $q$ and block size
as in the previous two examples. The weights are specified in the
following way:
\begin{itemize}
\item[-] NEW: $w_{o}=w_{s}=1/3$ and $w_{m,i}=1/12$, $i=1,\ldots,4$.
\item[-] PS: $w_{o}=w_{s}=0$ and $w_{m,i}=1/4$, $i=1,\ldots,4$.
\end{itemize}
%
%
\begin{table}
\caption{Distribution of patients among 160 strata}
\label{t120among160}
\begin{tabular*}{\textwidth}{@{\extracolsep{\fill}}lccccc@{}}
\hline
\textbf{\# of pts within stratum} &\textbf{0} &\textbf{1} & \textbf
{2} &\textbf{3} &\multicolumn{1}{c@{}}{\textbf{4 and more}}\\
\hline
\# of strata &95.4\phantom{\%} &38.8\phantom{\%} & 12.7 &5.6\phantom
{\%} &7.6\phantom{\%} \\
Proportion &59.6\% &24.3\% &\phantom{000}7.9\% &3.5\% &4.7\% \\
\hline
\end{tabular*}
\end{table}

Table~\ref{t120among160} shows the distribution of 120 patients among
160 strata. In this case 24.3\% of the strata have 1 patient; 11.4\%
contain 2 or 3 patients. If stratified randomization is employed, then
the patients in the above 24.3\% strata has to be randomized by equal
probabilities. Moreover, the incomplete blocks in strata with~2 or 3
patients also pose a high risk of large overall imbalance.

%
%
\begin{table}[b]
\caption{Comparison of absolute overall imbalance $|D_{n}|$}
\label{tsimoverall}
\begin{tabular*}{200pt}{@{\extracolsep{\fill}}ld{2.2}d{2.2}d{2.2}@{}}
\hline
&\multicolumn{1}{c}{\textbf{STR-PB}} &\multicolumn{1}{c}{\textbf
{PS}}&\multicolumn{1}{c@{}}{\textbf{NEW}}\\
\hline
Mean &6.70 &0.91 &0.63 \\
Median &6 &0 &0 \\
95\% quan &16 &2 &2 \\
\hline
\end{tabular*}
\end{table}

The mean absolute imbalances at the three levels are compared, as shown
in Tables~\ref{tsimoverall},~\ref{tsimmarginal} and~\ref{tsimstratum}.
Table~\ref{tsimoverall} shows the result for the overall imbalance
and lists the mean, median and 95\% quantile of $|D_{120}|$. It is
seen that NEW has mean, median and $95\%$ quantile of 0.63, 0 and 2,
respectively, whereas PS has slightly higher values. The three
quantities are extremely high under STR-PB, which are not recommended
for this case.

%
\begin{table}
\caption{Comparison of mean absolute marginal imbalances $E|D_{n}(i;
k_i)|$}
\label{tsimmarginal}
\begin{tabular*}{\textwidth}{@{\extracolsep{\fill}}lcccc@{}}
\hline
& &\textbf{STR-PB} &\textbf{PS}&\textbf{NEW} \\
\hline
Gender &male &5.52 &1.10 &1.59 \\
&female &3.86 &1.06 &1.55 \\[3pt]
Age &$<60$ &4.84 &1.08 &1.57 \\
&$\geq60$ &4.40 &1.11 &1.23 \\[3pt]
Disease &moderate &5.01 &1.10 &1.56 \\
&severe &4.35 &1.18 &1.52 \\[3pt]
20 sites &2 small &1.45 &0.94 &1.02 \\
&16 median &1.44 &1.21 &1.32 \\
&2 large &1.47 &1.33 &1.52 \\
\hline
\end{tabular*}
\end{table}

Table~\ref{tsimmarginal} gives the mean absolute marginal imbalances.
For the covariates of gender, age and disease, the table explicitly
lists the mean values on these 6 margins, as each of them only has two
levels. For example, over the 1000 simulations, on average the absolute
differences of patients in the two treatment groups within all male are
5.52, 1.10 and 1.59 under STR-PB, PS and NEW, respectively. Therefore,
in this respect PS has the best performance; NEW is slightly worse, but
still tolerable; STR-PB is the worst, since its mean is as high as
5.52. Similar conclusions can be reached for the other 5 margins.
Moreover, for the margins relating to ``site,'' since there are a total
of 20 margins, we are unable to show the result on each margin due to
the space limit. Hence, these 20 margins are further categorized into
three groups of small, median and large sizes, and the mean values in
the table are further averaged over the margins within the groups. For
example, 1.32 is the mean absolute imbalance over the 16 median-sized
sites as well as over the 1000 simulations. In terms of imbalances on
margins defined by site, PS is still the best, and STR-PB has similar
performance to NEW. This is because each margin of site contains only 8
strata, hence the ``accumulating effect'' of within-stratum imbalances
under STR-PB is not as strong.

%
\begin{table}[b]
\caption{\mbox{\hspace*{-3pt}Comparison of absolute within-stratum imbalances
$|D_{n}(k_1,\ldots,k_I)|$: Distribution and mean}}
\label{tsimstratum}
\begin{tabular*}{\textwidth}{@{\extracolsep{\fill}}lcccc@{}}
\hline
\textbf{\# of pts' within strt.} & $\bolds{|D_{n}(k_1,\ldots,k_I)|}$
&\textbf{STR-PB} &\textbf{PS}&\textbf{NEW}\\
\hline
2 &prob($=0$) &0.68 &0.57 &0.69\\
&prob($=2$) &0.32 &0.43 &0.31 \\
& mean &\textbf{0.64} &\textbf{0.86} &\textbf{0.62} \\[3pt]
3 &prob($=1$) &1.00 &0.85 &0.94 \\
&prob($=3$) &0.00 &0.15 &0.06 \\
&mean &\textbf{1.00} &\textbf{1.30} &\textbf{1.12} \\
\hline
\end{tabular*}
\end{table}

Table~\ref{tsimstratum} displays the distribution and absolute mean
of within-stratum imbalances for strata with 2 or 3 patients. For
example, of all the strata which contain 2 patients, the absolute
difference is either 0 or 2, and the distribution is 0.69 to 0 and 0.31
to 2 under NEW, leading to an average of 0.62. According to this
criterion, NEW has the lowest mean, STR-PB has a~slightly larger value
and PS has mean as large as 0.86. For strata containing 3 patients,
since the block size is 4 for STR-PB, it is impossible to get an
absolute value of 3. Hence, the mean absolute imbalance is 1, the
minimum among the three methods.

In summary, our new method maintains good balance from all three
perspectives and should be favored. We also performed the simulations
under other parameter values. Some of them include: (1) changing the
weights $w_{o}$, $w_{s}$, and $w_{m,i}$, as well as the block size; (2)
$2\times100$ strata, representing few covariates but many levels at
least for one covariate; (3) $3\times4 \times5 \times6$ strata,
representing a few covariates and a few levels for each. In all the
above settings, our new procedure shows advantages over the other two methods.

\section{Conclusion}
\label{sconclusion}

In this paper we propose a new covariate-adaptive design that minimizes
a weighted average of
three types of imbalances (within-stratum, within-covariate-margin and
overall). Simulation results show that the proposed method provides
better allocation balance from different perspectives, while stratified
randomization and Pocock and Simon's (\citeyear{PocSim75}) marginal
method have large imbalances either as a whole, or within-stratum.

The new procedure can also be generalized in several ways. In this
paper we only considered balanced allocation (1:1), whereas in some
problems unequal ratios [\citet{HuRos06}] are also desired. For
example, if the two groups are an innovation versus a placebo, and a
pilot study has shown some effect of the innovation, then it is more
ethical to assign more patients to the innovation. If one treatment is
much more costly than the other, then assigning more patients to the
latter would be more economical. Sometimes, the randomization has to be
adapted to covariates as well as responses. \citet{Zhaetal07}
proposed ``covariate-adjusted response-adaptive randomization,'' whose
allocation ratio depends on both covariate profiles and responses of
patients. One may modify our proposed procedure to accommodate these
situations. On the other hand, some trials (e.g., some Phase~II trials)
involve the comparison of more than two treatments [\citet
{PocSim75}, \citet{HuRos06}, etc.]. We can generalize the proposed
procedure to clinical trials for comparing three or more treatments. We
leave these as future research topics.

For Efron's (\citeyear{Efr71}) biased coin design (without involving
covariates), it is well known that the imbalance is a positive
recurrent Markov chain. \citet{MarRos10} studied some exact
properties of Efron's (\citeyear{Efr71}) biased coin design.
However, to our best knowledge, there is no theoretical result about
the imbalance of covariate-adaptive randomization in literature, due to
the complex of the problem and the lack of technical tools. In this
paper, we introduced the technique of ``drift conditions'' in Markov
chains to study the theoretical properties of covariate-adaptive randomization.
This technique could provide a possible way of studying the properties
of general covariate-adaptive designs as well as covariate-adjusted
response-adaptive designs.

The inference under covariate-adaptive randomization is also an
important issue. By simulation studies, several authors have raised
concerns about the conservativeness of the unadjusted analysis (such as
two-sample $t$-test) under covariate-adaptive randomization and suggested
that all covariates that are used in the randomization should be
included in the analysis [\citet{Bir85}, \citet{For87},
etc.]. \citet{ShaYuZho10} studied the theoretical relationship
between different randomization designs and different inference
methods. To make the problem more tractable, the authors focused on a
simple homogeneous linear model. They found that if the underlying
response-covariate model can be correctly specified, then the usual
regression analysis is valid and has the highest power as compared to
other types of analysis, no matter what randomization is employed.
These results also apply to the proposed randomization procedure in
this paper.

If the model specification is not feasible and only a two-sample $t$-test
can be used, then the test under stratified randomization tends to have
a conservative type~I error rate due to the overestimation of
$\operatorname{\mathsf{Var}}
(\bar
{Y}_1-\bar{Y}_2)$. \citet{ShaYuZho10} used a bootstrap method to
correct the variance estimation. The resulting bootstrap $t$-test
restores the type I error rate, and is more powerful than the
traditional $t$-test under simple randomization. Similar bootstrap
adjustment can be used as an inference method for the new randomization
procedure. We leave this as a future research project.

\section{Sketch of proofs}
\label{ssketch}
\mbox{}
\begin{pf*}{Proof of Theorem~\ref{theorem1}}
With $2\times2$ strata, the within-stratum imbalances $\mathbf{D}_{n}$
and the multinomial probabilities $\mathbf{p}$ are both matrices of
$2\times2$. Let $\tilde{\mathbf{D}}_{n}=(D_{n,1},D_{n,2},D_{n,3},D_{n,4})
:=(D_n(1,1),D_n(1,2),D_n(2,1),D_n(2,2))$, that is, $\tilde{\mathbf
{D}}_{n}$ is simply the vector form of $\mathbf{D}_{n}$. $\tilde
{\mathbf
{p}}=(p_1,\ldots,p_{4})$ can be defined in the same way. By the above
notation, any stratum can be represented by the $2$-index form
$(k_{1}, k_{2})$, or the single-index form $(r)$ ($1\leq r \leq4$).
The quantity $\delta_{n-1}(k_1^*,k_2^*)$ in (\ref{simplification}) then
reduces to
%
%
\begin{eqnarray}\label{DtoDelta}
&&\delta_{n-1} \bigl(k_1^*,k_2^* \bigr)\nonumber\hspace*{-25pt}
\\[-1pt]
&&\quad=(w_{o}+w_{m,1}+w_{m,2}+w_{o})D_{n-1}
\bigl(k_1^*,k_2^* \bigr)+(w_{o}+w_{m,1})D_{n-1}
\bigl(k_1^*,k_2 \bigr)
\nonumber\hspace*{-25pt}
\\[-8.5pt]
\\[-8.5pt]
\nonumber
&&\qquad{}+(w_{o}+w_{m,2})D_{n-1} \bigl(k_1,k_2^*
\bigr)+w_{o}D_{n-1}(k_1,k_2)\hspace*{-25pt}
\\[-1pt]
&&\quad=u_{1}D_{n-1} \bigl(k_1^*,k_2^*
\bigr)+u_{2}D_{n-1} \bigl(k_1^*,k_2
\bigr)+u_{3}D_{n-1} \bigl(k_1,k_2^*
\bigr)+u_{4}D_{n-1}(k_1,k_2),
\nonumber\hspace*{-25pt}
\end{eqnarray}
where $k_1\neq k_1^*$, $k_2\neq k_2^*$ and $u_{1}=1$. Let $\tilde
{\bolds{\delta}}_{n}=(\delta_{n,1},\delta_{n,2},\delta_{n,3},\delta_{n,4})
:=(\delta_n(1,1),\break \delta_n(1,2),$ $\delta_n(2,1),\delta_n(2,2))$. Then,
according to (\ref{DtoDelta}), $\tilde{\mathbf{D}}_{n}$ and $\tilde
{\bolds{\delta}}_{n}$ are linked~by
\begin{equation}
\tilde{\bolds{\delta}}_{n}=\tilde{\mathbf{D}}_{n}\pmatrix{
u_{1} & u_{2} & u_{3} & u_{4}
\vspace*{2pt}
\cr
u_{2} & u_{1} & u_{4} &
u_{3} \vspace*{2pt}
\cr
u_{3} & u_{4} &
u_{1} & u_{2} \vspace*{2pt}
\cr
u_{4} &
u_{3} & u_{2} & u_{1} \vspace*{2pt}
\cr
}:=
\tilde{\mathbf{D}}_{n}\mathbf{U}.\label{matrixDtoDelta}\vadjust{\goodbreak}
\end{equation}
For any $\tilde{\mathbf{D}}_{n}\in\mathbb{Z}^{4}$, we define a test function
\[
V(\tilde{\mathbf{D}}_{n})=\sum_{r=1}^{4}
\frac{[D_{n,r}]^{2}}{p_{r}},
\]
that is, the sum of squared within-stratum differences adjusted for the
corresponding multinomial
probabilities. The test set $\mathcal{C} $ is defined as $\mathcal
{C}=\{
\tilde{\mathbf{D}}_{n}\dvtx\max_{r}\Vert\tilde{D}_{n,r}\Vert
\leq K\}$ ($K>0$ is
to be determined). $V$ and $\mathcal{C} $ are the key elements in
proving positive recurrence, according to the drift conditions~(\ref{drift1}) and (\ref{drift2}).

For the ease of representation, in the rest of the proof we will simply
use the notation $\mathbf{D}$ and $\bolds{\delta}$ for $\tilde
{\mathbf{D}}_{n}$ and $ \tilde{\bolds{\delta}}_{n}$, respectively,
unless specified otherwise. Under the new allocation rule, it can be
derived that the one-step movement $\Delta V(\mathbf{D})$, defined in
(\ref{deltavd}), has the form
\[
\Delta V (\mathbf{D})=2(q-p)\sum_{r=1}^{4}\operatorname{sgn}(
\delta_{r})D_{r}+4,
\]
where $D_{r}$ and $\delta_{r}$ are the $r$th element of vectors
$\mathbf
{D}$ and $\bolds{\delta}$, respectively, and $\operatorname{sgn}(x)=1$, $-1$, $0$
for $x>0$, $<0$ or $=0$. For derivation of $\Delta V (\mathbf{D})$, see
Section~1 of the supplemental article [\citet{HuHu}].

We need to show that $\Delta V (\mathbf{D})$ satisfies drift conditions
(\ref{drift1}) and (\ref{drift2}). In fact, since the test set
$\mathcal
{C}$ is bounded, (\ref{drift2}) is trivially true. Since $q-p<0$,
(\ref
{drift1}) is equivalent to finding $M'_{1}>2/(p-q)$ such that
%
%
\begin{equation}
\Delta W(\mathbf{D}):=\sum_{r=1}^{4}\operatorname{sgn}(
\delta_{r})D_{r}> M'_{1}\qquad  \mbox{for } \mathbf{D}\notin\mathcal{C}.\label{drift4}
\end{equation}
Intuitively, when $u_{2}$, $u_{3}$ and $u_{4}$ are small, $\delta_{r}$
is expected to be close to $D_{r}$ so that they have the same sign.
Thus, a larger proportion of the strata have
$\operatorname{sgn}(\delta_{r})D_{r}=\operatorname{sgn}(D_{r})D_{r}=|D_{r}|$ and $\Delta W(\mathbf
{D})$ tends to be positive. In the trivial case that
$u_{2}=u_{3}=u_{4}=0$, that is, $\mathbf{D}=\bolds{\delta}$, we have
$\Delta W(\mathbf{D})=\sum_{r=1}^{4}|D_{r}|>K$, so (\ref{drift4}) holds
by letting $M'_{1}=K=2.1/(p-q)$. Therefore, in the following we can
assume that $\max\{u_{2},u_{3},u_{4}\}>0$.

For any $\mathbf{D}\in\mathcal{C}^c\subset\mathbb{Z}^{4}$, call the
pair of $(D_{r},\delta_{r})$ a ``\textit{match}'' if
$\delta_{r}\neq0$ and $\delta_{r}D_{r}\geq0$.
Hence, for a match $\operatorname{sgn}(\delta_{r})D_{r}=|D_{r}|$.
Furthermore, define $M(\mathbf{D},\bolds{\delta})$ as the number
of matches
in $(D_{r},\delta_{r})$'s, $r=1,\ldots,4$. Obviously,
$0\leq M(\mathbf{D},\bolds{\delta})\leq4$. It can be shown that
$M(\mathbf{D},\bolds{\delta})=0$ is impossible for $\mathbf{D}\in\mathcal{C}^c$.
Therefore, for $M(\mathbf{D},\bolds{\delta})=i$,
$i=1,2,3,4$, if we can find $d_{i}>0$, such that $\Delta W(\mathbf
{D})>d_{i}K$, then (\ref{drift4}) is true by letting $M_1'=K\min_{i}d_{i}$ and
$K=2.1/[(p-q)\min_{i}d_{i}]$.

In fact, finding $d_4$ for $M(\mathbf{D},\bolds{\delta})=4$ is
quite trivial ($d_{4}=1$). We will show how to find $d_3$ for
$M(\mathbf
{D},\bolds{\delta})=3$ below. When $M(\mathbf{D},\bolds
{\delta
})=3$, we know that $a_{1}=\max\{u_{2},u_{3},u_{4}\}\neq0$ and
$a_{2}=\min\{1-u_{2},1-u_{3},1-u_{4}\}\neq0$ (since \mbox{$w_{s}\neq0$}).
Without loss of generality assume $D_{1}$ and $\delta_{1}$ do not
match, which means $\delta_{1}D_{1}\leq0$. Thus\vadjust{\goodbreak} $|\delta_{1}-D_{1}|\geq
|D_{1}|$. By (\ref{matrixDtoDelta}), $\delta_{1}-D_{1}=u_{2}D_{2}+u_{3}D_{3}+u_{4}D_{4}$, which implies
$|u_{2}D_{2}+u_{3}D_{3}+u_{4}D_{4}|\geq|D_{1}|$.
Then
\begin{eqnarray*}
\Delta W(\mathbf{D})&\geq&-|D_{1}|+|D_{2}|+|D_{3}|+|D_{4}|
\\[-2pt]
&\geq&-\bigl(u_{2}|D_{2}|+u_{3}|D_{3}|+u_{4}|D_{4}|\bigr)+|D_{2}|+|D_{3}|+|D_{4}|
\\[-2pt]
&\geq& a_{2}\bigl(|D_{2}|+|D_{3}|+|D_{4}|\bigr)
\\[-2pt]
&\geq& a_{2} \bigl[(1/2) \bigl(|D_{2}|+|D_{3}|+|D_{4}|\bigr)+(1/2)a_{1}^{-1}|D_{1}|
\bigr]
\\[-2pt]
&\geq&(a_{2}/2)\min \bigl\{1, a_{1}^{-1} \bigr\}
\cdot\max\bigl\{ |D_{1}|,|D_{2}|,|D_{3}|,|D_{4}|
\bigr\}
\\[-2pt]
&> &(a_{2}/2)\min \bigl\{1, a_{1}^{-1} \bigr\}
\cdot K := d_{3}K.
\end{eqnarray*}
The ways of finding $d_2$ and $d_1$ for $M(\mathbf{D},\bolds
{\delta
})=2$ and $1$ are similar, but require more work. In particular,
condition \textup{(B)} in Theorem (\ref{theorem1}) is needed to verify the
case of $M(\mathbf{D},\bolds{\delta})=1$. In Section~2 of the
supplemental article [\citet{HuHu}], we show how to find $d_i$ for
$i=4,3,2,1$ and explain why $M(\mathbf{D},\bolds{\delta})\neq0$.

Corollary~\ref{corollary1} is obtained by solving the linear equation
for $\mathbf{x}$ in Theorem~\ref{theorem1} under the assumption that
$w_{m,1}=w_{m,2}$ and then substituting the solution in
$|x_{1}|+|x_{2}|+|x_{3}|<1$. For detailed proof of Corollary~\ref
{corollary1}, see Section~3 of the supplemental article [\citet{HuHu}].
\end{pf*}

\begin{pf*}{Proof of Theorem~\ref{theorem2}} The main steps are
similar to those in Theorem~\ref{theorem1}. Let $\tilde{\mathbf
{D}}_{n}=(D_{n,1},\ldots, D_{n,m})$ be the vector version of $\mathbf
{D}_{n}=(D_{n}(k_{1},\ldots,\break k_{I}))$, where the $m$ strata can be
arbitrarily ordered and indexed by $1,\ldots, m$.
Similary, let $\tilde{\bolds{\delta}}_{n}$ and ~$\tilde{\mathbf
{p}}$ be the vector forms of array $(\delta_n(k_1,\ldots,k_I))$ and
array $(p(k_1,\ldots,k_I))$, respectively, using the same order as in
$(D_{n}(1),\ldots, D_{n}(m))$. By the above notation, any stratum can
be represented by the $I$-index form $(k_{1},\ldots, k_{I})$, or the
single-index form $(r)$ ($1\leq r \leq m$). As in the $2\times2$ case,
let $\tilde{\bolds{\delta}}_{n}:=\tilde{\mathbf
{D}}_{n}\mathbf
{U}$. Then by the definition of $\tilde{\bolds{\delta}}_{n}$ as
well as the description of weights before Theorem~\ref{theorem2}, for
any two strata $(r)=(k_{1}^{*},\ldots, k_{I}^{*})$ and
$(s)=(k_{1},\ldots, k_{I})$, the element $u_{rs}$ in the matrix of
$\mathbf{U}$ is determined as follows:
for any fixed $i_{1}< i_{2}< \cdots< i_{l}$ ($\{i_{1},\ldots,i_{l}\}
\in
\mathbb{I}$), if $k_{i_{t}}\neq k_{i_{t}}^{*}$ and $k_{j}=k_{j}^{*}$
for $j\in\mathbb{I}$, $j\neq i_{t}$ and $1\leq t \leq l$, then
\[
u_{rs}=w_{o}+\sum_{j\neq i_{t}, 1\leq t \leq l}w_{m,j}.
\]
So $u_{rs}=u_{sr}$, and for any $r$, $\sum_{s=1,\ldots,m, s\neq
r}u_{rs}=u^{*}$, as defined in Theorem~\ref{theorem2}.

The test function $V$ and the test set $\mathcal{C}$ are still defined
as before, except that in this case the dimension of $\tilde{\mathbf
{D}}_{n}$ is $m$ instead of $4$. Use the simple notation $\mathbf{D}$
and $\bolds{\delta}$ for $\tilde{\mathbf{D}}_{n}$ and $ \tilde
{\bolds{\delta}}_{n}$, respectively. In the same manner, to verify
the drift conditions it is equivalent to find $M'_1>2/(p-q)$ such that
%
%
\begin{equation}
\Delta W(\mathbf{D}):=\sum_{r=1}^{m}\operatorname{sgn}(
\delta_{r})D_{r}> M'_{1} \qquad\mbox{for } \mathbf{D}\notin\mathcal{C}.\label{drift5}
\end{equation}
For any fixed $\mathbf{D}\in\mathcal{C}^c\subset\mathbb{Z}^{m}$,
suppose for $(D_{r},\delta_{r})$'s, $r=1,\ldots,m$, there are $m_{0}$
mismatched pairs. Without loss of generality assume that the mismatched
pairs occur in\vadjust{\goodbreak} the $1$st, $2$nd, $\ldots$ and the $m_{0}$th strata. By
the definition of a mismatched pair, $D_{r}$ and $\delta_{r}=\sum_{s=1}^{r-1}u_{rs}D_{s}+D_{r}+\sum_{s=l+1}^{m}u_{rs}D_{s}$ have
different signs, $r=1,\ldots, m_{0}$. Therefore,
%
%
\begin{eqnarray}\label{drs}
|D_{r}|&\leq&\Biggl|\sum_{s=1}^{r-1}u_{rs}D_{s}+
\sum_{s=l+1}^{m}u_{rs}D_{s}\Biggr|
\nonumber
\\[-9pt]
\\[-9pt]
\nonumber
&\leq&\sum_{s=1,\ldots, m_{0} ,s\neq r}u_{rs}|D_{s}|+
\sum_{s=m_{0}+1}^{m}u_{rs}|D_{s}|.
\end{eqnarray}
First, we notice that $m_{0}\neq m$; otherwise, by summing (\ref{drs})
over $r=1$ to~$m$, we have $\sum_{r=1}^{m}|D_r|\leq u^*\sum_{r=1}^{m}|D_r|$ which
is impossible for $\mathbf{D}\in\mathcal{C}^c$
and $u^*<1/2$. Second, suppose $m_{0}\!\neq\!0$. By summing (\ref{drs})
over $r\!=\!1$ to $r\!=\!m_{0}$, we have
\[
\sum_{r=1}^{m_{0}}\biggl(1-\sum
_{s=1,\ldots, m_{0},s\neq
r}u_{rs}\biggr)|D_{r}|\leq \sum
_{r=m_{0}+1}^{m}\Biggl(\sum
_{s=1}^{m_{0}}u_{rs}\Biggr)|D_{r}|.
\]
Combined with the fact that $1-u^{*}\leq(1-\sum_{s=1,\ldots,
m_{0},k\neq r}u_{rs})$ for $r=1,\ldots,m_{0}$
and $(\sum_{s=1}^{m_{0}}u_{rs})\leq u^{*}$ for $r=m_0+1,\ldots,m$, it
is seen that
\[
\sum_{r=1}^{m_{0}}\bigl(1-u^{*}
\bigr)|D_{r}|\leq\sum_{r=m_{0}+1}^{m}u^{*}|D_{r}|.
\]
Then
\begin{eqnarray*}
\Delta W(\mathbf{D})&\geq&-|D_{1}|-|D_{2}|-\cdots
-|D_{m_{0}}|+|D_{m_{0}+1}|+|D_{m_{0}+2}|+\cdots+|D_{m}|
\\
&\geq&-\frac{u^{*}}{1-u^{*}}\sum_{r=m_{0}+1}^{m}|D_{r}|+
\sum_{r=m_{0}+1}^{m}|D_{r}|= \biggl(1-
\frac{u^{*}}{1-u^{*}} \biggr)\sum_{r=m_{0}+1}^{m}|D_{r}|.
\end{eqnarray*}
Since $0\leq u^{*}<1/2$, we have $0<1-\frac{u^{*}}{1-u^{*}}\leq1$.
Hence, the above inequality is also true for $m_{0}=0$. If $\max_{m_{0}+1\leq r \leq m}|D_{r}|>K$,
then $\Delta W(\mathbf{D})>
(1-\frac
{u^{*}}{1-u^{*}})K$;
otherwise $\max_{1\leq r\leq m_{0}}|D_{r}|>K$ and
$\sum_{r=m_{0}+1}^{m}|D_{r}|\geq\frac{1-u^{*}}{u^{*}}\sum_{r=1}^{m_{0}}|D_{r}|$,
which means $\Delta W(\mathbf{D})> (1-\frac
{u^{*}}{1-u^{*}})\frac{1-u^{*}}{u^{*}}K>(1-\frac{u^{*}}{1-u^{*}})K$.
Thus, if we define $M'_{2}=(1-\frac{u^{*}}{1-u^{*}})K$ and $K=\frac
{2.1}{p-q}(1-\frac{u^{*}}{1-u^{*}})^{-1}$, then
\[
\Delta W(\mathbf{D})>M'_{2}>2/(p-q).
\]
\upqed\end{pf*}

\section*{Acknowledgments}
Special thanks go to anonymous referees, the Associate
Editor and the Editor for the constructive comments, which led to a~much improved version of the
paper.

\begin{supplement}[id=suppA]
\stitle{Additional proofs}
\slink[doi]{10.1214/12-AOS983SUPP} 
\sdatatype{.pdf}
\sfilename{aos983\_supp.pdf}
\sdescription{We provide additional proofs that are omitted in
Section~\ref{ssketch}. They include:
(1) derivation of $\Delta V(\mathbf{D})$; (2) the appropriate choice
of $d_i$
when $M(\mathbf{D},\bolds{\delta})=i$ ($i=4,3,2,1$); (3) proof of
Corollary~\ref{corollary1}.\vadjust{\goodbreak}}
\end{supplement}

%

%

\printaddresses

\end{document}